\documentclass[reqno,11pt]{amsart}
\usepackage{amsthm,amssymb,amsmath}
\usepackage{geometry}
\usepackage{xcolor}
\usepackage{hyperref}
\hypersetup{
	colorlinks=true,
	linkcolor=blue,
	filecolor=magenta,
	urlcolor=cyan,
	citecolor=teal,
}
\allowdisplaybreaks

\newtheorem{lem}{Lemma}[section]

\newtheorem{pro}[lem]{Proposition}
\newtheorem{thm}[lem]{Theorem}

\numberwithin{equation}{section}

\begin{document}
	\title[Metric properties of continued fractions]{Metric properties of continued fractions with large prime partial quotients}

	\author[W. Cheng]{Wanjin Cheng}
	\address[Wanjin Cheng]{School of Mathematics, South China University of Technology, Guangzhou, 510640, China}
	\email{chengwj0227@163.com}
	\author[W. Wu]{Wen Wu$^{\star}$}
	\address[Wen Wu]{School of Mathematics, South China University of Technology, Guangzhou, 510640, China}
	\email[corresponding author]{wuwen@scut.edu.cn}
	
	
	\begin{abstract}
		Let $x \in [0,1)$ with continued fraction expansion $[a_1(x),a_2(x),\dots]$, and let $\phi:\mathbb{N}\to\mathbb{R}^+$ be a non-decreasing function. We consider the numbers whose continued fraction expansions contain at least two partial quotients that are simultaneously large and prime, that is
		\[
		E'(\phi):=\Big\{x\in[0,1): \exists\, 1\leq k\neq l\leq n, \  a'_{k}(x),\ a'_{l}(x)\geq\phi(n)
		\ \text{for i.m. } n\in\mathbb{N}\Big\},
		\]
		where $a'_i(x)$ denotes $a_i(x)$ if $a_i(x)$ is prime and $0$ otherwise. We establish a zero-one law for the Lebesgue measure of $E'(\phi)$ and determine its Hausdorff dimension.
	\end{abstract}
	
	\maketitle

	\section{Introduction}
	Continued fractions are among the most powerful tools for studying Diophantine approximation, as their partial quotients directly characterize how well a real number can be approximated by rationals.
	Every $x\in[0,1)$ admits a finite or infinite continued fraction expansion of the form
	\begin{equation*}%\label{zhanshi}
		\begin{aligned}
			x &=\frac{1}{a_{1}(x)+\dfrac{1}{a_{2}(x)+\dfrac{1}{a_{3}(x)
						+\ddots}}}
			{\ =:\ }[a_1(x),a_2(x),a_3(x),\dots]
		\end{aligned}
	\end{equation*}
	where $a_1(x)=\lfloor1/x\rfloor$ and $a_{k}(x)=a_1(T^{k-1}(x))$ for $k\geq2$; these are called the partial quotients of $x$, and $T$ denotes the Gauss map (i.e., $T(x) = \{1/x\}$ for $x \neq 0$, $T(0)=0$).
	For $n\in\mathbb{N}$, the finite truncation
	$[a_{1}(x), a_{2}(x), \dots, a_{n}(x)]=:p_{n}(x)/q_{n}(x)$ is called the $n$-th convergent of $x=[a_1(x),a_2(x),a_3(x),\dots]$.
	It is well-known that for all $n\in \mathbb{N}$,
	$$\frac{1}{(a_{n+1}(x)+2)q_n^2}<\left|x-\frac{p_n}{q_n}\right|<\frac{1}{a_{n+1}(x)q_n^2},$$ see for example \cite{Khinchin}.
	This shows that the growth rate of partial quotients further quantifies the quality of approximation: if $a_k(x)$ is bounded for all $k$, then $x$ cannot be approximated ``too well''; conversely, unbounded partial quotients imply that $x$ is ``very well approximable''.
	
	Borel \cite{Borel} and Bernstein \cite{Bernstein} investigated the set of point with large partial quotients, defined as follows:
	\begin{equation*}
		F(\phi) := \big\{ x \in [0, 1) \colon a_{n}(x) \ge \phi(n) \ \text{for infinitely many } n \in \mathbb{N} \big\},
	\end{equation*}
	where $\phi\colon \mathbb{N} \to (0, \infty)$ is a positive function. They established a zero-one law for the Lebesgue measure of $F(\phi)$. In the case where $F(\phi)$ has zero Lebesgue measure, its Hausdorff dimension was subsequently studied by Good \cite{Good}, {\L}uczak \cite{Luczak}, and Wang--Wu \cite{Wang}.
	
	When $F(\phi)$ has full Lebesgue measure, it follows that for almost every $x\in [0,1)$, there exists a subsequence of the partial quotients of $x$ that tends to infinity at the rate of $\phi(n)$. In general, for Lebesgue almost all $x\in [0,1)$, the Ces\`aro mean of the partial quotients of $x$ diverges to infinity. That is, there no finite law of large numbers holds:
	\begin{equation*}
		\lim_{n \to \infty} \frac{1}{n}\sum^{n}_{i=1}a_i(x) = \infty
		\quad \text{for}\ \mathcal{L}\text{-a.e. } x\in(0,1).
	\end{equation*}
	Furthermore, Philipp \cite{Philipp} showed that no normalization yields a positive finite limit.
	However, Diamond and Vaaler \cite{Diamond} showed that by removing the maximal term, a precise asymptotic law emerges: for Lebesgue almost all $x \in [0,1)$,
	\begin{equation}\label{limit2}
		\lim_{n \to \infty} \frac{1}{n\log n}\Bigl(\sum^{n}_{i=1}a_i(x)-\max_{1\le k\le n}a_k(x)\Bigr)
		= \frac{1}{\log 2}.
	\end{equation}
	A key ingredient in proving \eqref{limit2} is the rarity of two large partial quotients occurring close to each other. This leads to the study of the set
	\[
	E(\phi) := \big\{ x \in [0, 1) \colon \exists\, 1 \le k \ne l \le n,\
	a_{k}(x) \ge \phi(n),\ a_{l}(x) \ge \phi(n)\ \text{for i.m. } n \in \mathbb{N} \big\},
	\]
	whose Lebesgue measure and Hausdorff dimension were completely characterized by Tan et al.\ \cite{TTW} for non-decreasing functions $\phi$.

	Schindler and Zweim\"uller \cite{Schindler} considered variations of the above questions by restricting the partial quotients to be primes; that is, by replacing $a_n(x)$ with $a_n'(x)$, where 
	\[a'_n(x)=\begin{cases}a_n(x), & \text{if } a_n(x) \text{ is prime},\\ 0, & \text{	otherwise}.\end{cases}\] 
	They proved the zero-one law for the Lebesgue measure of the set  
	\[	F'(\phi) := \big\{ x \in [0,1) \colon a_n'(x) \geq \phi(n) \ \text{for i.m. } n \in \mathbb{N} \big\},\] 
	where $\phi:\mathbb{N}\to\mathbb{R}^+$ is a positive function; namely, 
	\begin{equation}
		\mathcal{L}(F'(\phi))= \begin{cases}
			0, & \text{ if }\sum_{n\geq 1}\frac{1}{\phi(n)\log\phi(n)}<\infty,\\
			1, & \text{ if }\sum_{n\geq 1}\frac{1}{\phi(n)\log\phi(n)}=\infty,
		\end{cases}
		\label{0-1-[20]}
	\end{equation}
	see \cite[Theorem 2.1]{Schindler}. In the same setting, they also proved an analogue of the limit theorem \cite[Theorem 2.2]{Schindler} like \eqref{limit2} as follows: for Lebesgue almost all $x \in [0,1)$,
	\begin{equation}\label{limit2-prime}
		\lim_{n \to \infty} \frac{1}{n\log n}\Bigl(\sum^{n}_{i=1}a'_i(x)-\max_{1\le k\le n}a'_k(x)\Bigr)
		= \frac{1}{\log 2}.
	\end{equation}
	The rarity of two large prime partial quotients occurring close to each other plays an important role in proving \eqref{limit2-prime}.  
	This motivates the study of the ``size'' of the set 
	\[E'(\phi):=\big\{x\in[0, 1)\colon~\exists\ 1 \leq{k\neq l}\leq n,\ a_{k}'(x)\geq\phi(n),\ a_{l}'(x)\geq\phi(n) ~ \text{for~i.m.}~ n\in\mathbb{N} \big\}.\]

	Our first result is the zero-one law for the Lebesgue measure of $E'(\phi)$.
	\begin{thm}\label{t1} %Let $\phi$ be a positive function defined on natural numbers $\mathbb{N}$.
		Let $\phi:\mathbb{N}\rightarrow \mathbb{R}^+$ be a non-decreasing function. Then
		\begin{equation*}
			\mathcal{L}(E'(\phi))=\left\{
			\begin{aligned}
				0, &  \quad\quad  \text{if}\quad \sum_{n=1}^\infty\frac{n}{\phi^2(n)\log^2\phi(n)}<\infty ,\\
				1, & \quad\quad   \text{if} \quad \sum_{n=1}^\infty\frac{n}{\phi^2(n)\log^2\phi(n)}=\infty .
			\end{aligned}\right.
		\end{equation*}
	\end{thm}
	Our second result determines the Hausdorff dimension of $E'(\phi)$ when $\mathcal{L}(E'(\phi))=0$. We show that although \( E'(\phi) \subsetneq E(\phi) \), the Hausdorff dimension of $E'(\phi)$ is the same as that of $E(\phi)$ obtained by Tan et al.\ \cite{TTW}.
	\begin{thm}\label{t2}
		Let $\phi:\mathbb{N}\rightarrow \mathbb{R}^+$ be a non-decreasing function. Define
		$$\log B=\liminf_{n\rightarrow \infty}\frac{\log\phi(n)}{n}, \quad \log b=\liminf_{n\rightarrow \infty}\frac{\log\log \phi(n)}{n}.$$
		We have
		\[
		\dim_{H}E'(\phi)=\begin{cases}
			\inf\{s\geq 0: P(T,f)\leq0\}, & \text{if }1\leq B<\infty;\\
			\frac{1}{b+1}, & \text{if }B=\infty,
		\end{cases}
		\]
		where $f(x)=-(3s-1)\log B-s\log|T'(x)|$ and \(P(T,f)\) denotes the pressure function defined in \eqref{PT1}.
	\end{thm}

	This paper is organized as follows. In Section 2, we present the necessary preliminaries, including the elementary properties of Hausdorff dimension, continued fractions, and pressure functions, as well as certain results on prime numbers. In Section 3, we establish the zero-one law for the Lebesgue measure of \( E'(\phi) \), which is Theorem \ref{t1}. In Section 4, we prove Theorem \ref{t2} for \( 1 < B < \infty \), and address the cases \(B=1\) and \(B=\infty\) in Section 5.

	\section{Preliminaries}
	We first introduce some notation.
	For two positive real-valued functions \( f \) and \( g \), we write \( f \lesssim g \) if there exists a constant \( C > 0 \) such that \( f(x) \leq C g(x) \) for all \( x \).
	Similarly, we write \( f \asymp g \) if both \( f \lesssim g \) and \( f \gtrsim g \) hold.
	
	\subsection{The distribution of prime numbers}We recall several classical results in analytic number theory concerning the distribution of prime numbers.
	%	\begin{thm}[Prime Number Theorem \cite{Montgomery}]
		%		Let $\pi(x)$ be the prime-counting function, i.e. the number of prime numbers not exceeding $x$,
		%		we have
		%		$$\pi(n)\sim \frac{n}{\log n}.$$
		%	\end{thm}
	
	By the Prime Number Theorem, we can estimate the number of primes contained in various intervals. 
	%The following proposition characterizes the number of prime numbers in both large and small intervals.
	\begin{pro}[\cite{Robert,Bordignon}]\label{p2.5}
		The following statements hold:
		\begin{enumerate}
			\item[\rm(1)]\rm{(Large interval)}
			For any $\gamma > 1$ and $n \in \mathbb{N}$, define $c_n(\gamma)$ by
			\[
			\#\big(\mathcal{P} \cap [\gamma^n, 2\gamma^n]\big)
			= c_n(\gamma) \, \frac{\gamma^n}{n \log \gamma}.
			\]
			Then $c_n(\gamma) \to 1$ as $n \to \infty$.
			
			\item[\rm(2)]\rm{(Small interval)}
			If $x > e^{20}$, then $\mathcal{P} \cap [0.999x, x) \neq \emptyset$.
		\end{enumerate}
	\end{pro}
	
	The following proposition was observed by Schindler and Zweim\"uller \cite{Schindler} and plays a crucial role in the proof of Theorem \ref{t2}.
	\begin{pro}[\cite{Schindler}]\label{p2.8}
		For sufficiently large $M \in \mathbb{N}$, we have
		$$\sum_{p\geq M , \ p\in \mathcal{P}}\frac{1}{p^2}\asymp(\frac{1}{M\log M}).$$
	\end{pro}
	
	\subsection{Continued fractions}
	%\cite{Khinchin,Hardy,Iosifescu} and references therein.
	For $n\in\mathbb{N}$, call
	$p_{n}(x)/q_{n}(x)=[a_{1}(x), a_{2}(x), \dots, a_{n}(x)]$ the  $n$-th convergent of the continued fraction of $x$.
	With the conventions  $p_{-1}(x)=1$, $q_{-1}(x)=0$, $p_{0}(x)=0$ and $q_{0}(x)=1$, we have
	\begin{equation*}
		p_{n+1}(x)=a_{n+1}(x)p_{n}(x)+p_{n-1}(x),\quad
		q_{n+1}(x)=a_{n+1}(x)q_{n}(x)+q_{n-1}(x).
	\end{equation*}
	
	For any $n\geq1$, and $\big(a_{1}, a_{2}, \dots, a_{n}\big)\in\mathbb{N}^{n}$, let $q_{n}(a_{1},a_{2},\dots,a_{n})$ be the denominator of the finite continued fraction $[a_{1},a_{2},\dots,a_{n}]$. If there is no confusion, we write $q_{n}$ instead of $q_{n}(a_{1},a_{2},\dots,a_{n})$ for simplicity. %$q_{n}\geq2^\frac{n-1}{2} \text{and}\ q_n\leq(a_n+1)q_{n-1}.$
	\begin{lem}[\cite{Khinchin}]{\label{eq:2}}
		For any $n\geq1$, $k\geq1$ and $\big(a_{1}, a_{2}, \dots, a_{n}\big)\in\mathbb{N}^{n}$, we have
		\begin{enumerate}
			\item$q_n\geq 2^{\frac{n-1}{2}}$.
			\item$1\leq \frac{q_{n+k}(a_1,\dots,a_n,\dots,a_{n+k})}{q_n(a_1,\dots,a_n)q_k(a_{n+1},\dots,a_{n+k})}\leq2.$
			%\item$\frac{a_k+1}{2}\leq \frac{q_n(a_1,\dots,a_n)}{q_{n-1}(a_1,\dots,a_{k-1},a_{k+1},\dots,a_n)}\leq a_{k+1}.$
		\end{enumerate}
	\end{lem}
	For any $n\geq1$ and $\mathbf{a}=(a_1,\dots,a_n)\in\mathbb{N}^n$, denote by $I_{n}(a_{1}, a_{2}, \dots, a_{n})$ (or $I_n(\mathbf{a})$ for simplicity) the $n$-th cylinder of continued fractions	$$I_n(\mathbf{a})=I_{n}(a_{1}, a_{2}, \dots, a_{n}):=
	\big\{x\in[0,1)\colon a_{1}(x)=a_{1}, a_{2}(x)=a_{2}, \dots, a_{n}(x)=a_{n}\big\}.$$
	The following lemmas describe the length and the Lebesgue measure of the cylinder \(I_n(\mathbf a)\).
	\begin{lem}[\cite{Khinchin}]\label{L2.9} For any $n\geq1$, and $\mathbf{a}=\big(a_{1}, a_{2}, \dots, a_{n}\big)\in\mathbb{N}^{n}$, we have
		%		\begin{equation*}
			%			I_{n}\big(a_{1}, a_{2}, \ldots, a_{n}\big)=\left\{
			%			\begin{aligned}
				%				&\left[\frac{p_{n}}{q_{n}},\frac{p_{n}+p_{n-1}}{q_{n}+q_{n-1}}\right),
				%				\quad\quad  \text{if $n$ is even} ,\\
				%				&\left(\frac{p_{n}+p_{n-1}}{q_{n}+q_{n-1}}, \frac{p_{n}}{q_{n}}\right],
				%				\quad\quad  \text{if $n$ is odd},
				%			\end{aligned}\right.
			%		\end{equation*}
		%		and its length satisfis%a direct calculation gives the length of $I\big(a_{1}, a_{2}, \ldots, a_{n}\big)_{N}$
		\begin{equation*}
			\frac{1}{2q_{n}^2}\leq\big|I_{n}(\mathbf{a})\big|=\bigg|\frac{p_{n}+p_{n-1}}{q_{n}+q_{n-1}}-
			\frac{p_{n}}{q_{n}}\bigg|=\frac{1}{q_{n}\big(q_{n}+q_{n-1}\big)}\leq\frac{1}{q_{n}^2}.
		\end{equation*}
	\end{lem}
	\begin{lem}[\cite{Khinchin}]\label{lemma2.11}Let $I_{n}(\mathbf{a})$ be an $n$-th cylinder with $\mathbf{a}\in\mathbb{N}^n$, and let $M\in\mathbb{N}$. Then
		$$\frac{\mathcal{L}(I_n(\mathbf{a}))}{3M^{2}}\leq \mathcal{L}\big(I
		_{n+1}(\mathbf{a}, M)\big) \leq\frac{2\mathcal{L}(I_n(\mathbf{a}))}{M^{2}}.$$
		% $$\frac{\mathcal{L}(I_n)}{3M^{2}}\leq \mathcal{L}\bigg(\bigcup_{a_{n+1}=M}I\big(a_{1}, \ldots, a_{n}, a_{n+1}\big)\bigg) \leq\frac{2\mathcal{L}(I_n)}{M^{2}}.$$
	\end{lem}
	
	\subsection{The measure of a limsup set.} Chung-Erd\"{o}s \cite{Chung} provided a lower bound for the measure of limsup sets.
	\begin{lem}[Chung-Erd\"{o}s Inequality \cite{Chung}]\label{L2.15}
		Let $(\Omega,\mathcal{A},\mu)$ be a finite measure space, and let $\{E_n\}_{n\geq 1}$ be a countable family of measurable sets in $\mathrm{\Omega}$ with $\sum_{n=1}^{\infty}\mu(E_n)=\infty$. Then
		$$\mu\Bigl(\limsup\limits_{n\rightarrow\infty}E_n\Bigr)\geq \limsup\limits_{n\rightarrow\infty}\frac{(\sum_{i=1}^{n}\mu(E_i))^2}{\sum_{i,j=1}^{n}\mu(E_i\cap E_j)}.$$
	\end{lem}
	To derive the full measure, we need the following lemma. %to apply the Chung--Erd\"{o}s lemma locally and then the following Inequality allows us to conclude full measure.
	\begin{lem}[{\cite[Lemma 3.3.1]{Dajani}}, \cite{Knopp}]\label{L2.16}
		Let $I\subset \mathbb{R}$ be a bounded interval. If $E\subset I$ is a Lebesgue measurable set and $\mathcal{C}$ is a class of subintervals of $I$ satisfying
		\begin{enumerate}
			\item every open subinterval of $I$ is at most a countable union of disjoint elements from $\mathcal{C}$,
			\item for any $C\in \mathcal{C}$, $\mathcal{L}(C\cap E)\geq \rho\mathcal{L}(C)$, where $\rho>0$ is a constant independent of $C$, 
		\end{enumerate}
		then $\mathcal{L}(E)=\mathcal{L}(I)$.
	\end{lem}
	
	\subsection{Pressure function.} We end this section with a definition of the pressure function. For more details, see \cite{Hanus, Mauldin1, Mauldin2}.
	
	For any non-empty subset $A \subset \mathbb{N}$, define
	\[
	X_{A} = \{ x \in [0,1) : a_n(x) \in A \ \text{for all}\ n \ge 1 \}.
	\]
	Then $(X_{A}, T)$ is a sub-system of $([0,1),T)$.
	
	Let $g: [0,1) \to \mathbb{R}$ be a real-valued function.
	The pressure function on $(X_A, T)$ is defined by
	\begin{equation}\label{PT1}
		P_{A}(T, g)
		= \lim_{n \to \infty}
		\frac{1}{n}
		\log
		\sum_{\mathbf{a} \in A^n}
		\exp\bigl(\sup_{x \in X_{A} \cap I_n(\mathbf{a})}
		S_n g(x)\bigr),
	\end{equation}
	where $S_n g(x) = \sum_{j=0}^{n-1} g(T^j x)$.
	In particular, we denote $P(T, g) := P_{\mathbb{N}}(T, g)$. The existence of the limit in \eqref{PT1} was proved in \cite{Li}.
	%The $n$-th variation of $\varphi$ is defined as
	%\[
	%\mathrm{Var}_n(\varphi) :=
	%\sup\{\, |\varphi(x) - \varphi(y)| : I_n(x) = I_n(y) \,\}.
	%\]
	
	%The following proposition guarantees the existence of the limit in the above definition of the pressure function.
	
	%\begin{pro}[\cite{Li}]\label{exsit}
	%Suppose $\varphi: [0,1) \to \mathbb{R}$ satisfies
	%$\mathrm{Var}_1(\varphi) < \infty$ and
	%$\lim_{n \to \infty} \mathrm{Var}_n(\varphi) = 0$.
	%Then the limit in \textup{\eqref{PT1}} exists.
	%In particular, the pressure function is well-defined, and the value of $P_A(T, \varphi)$ remains unchanged even if the supremum in \textup{\eqref{PT1}} is omitted.
	%\end{pro}
	
	%Next, we present a continuity property of the pressure function in the continued fraction setting:
	%the pressure for the full system $([0,1), T)$ can be approximated by the pressures of its subsystems $(X_A, T)$ (see \cite{Li} for further details).
	
	%\begin{pro}\label{cor2.17}
	%Under the assumptions above, we have
	%\[
	%P_{\mathbb{N}}(T, \varphi)
	%=
	%\sup\big\{ P_A(T, \varphi) : A \text{ is a finite subset of } \mathbb{N} \big\}.
	%\]
	%\end{pro}

	\section{Proof of Theorem \ref{t1}}
	
	%	Before proving Theorem~\ref{t1}, we make a simple observation concerning the structure of the set \( E'(\phi) \). 
	Recall that
	\[E'(\phi)
	:= \big\{ x \in [0,1) : \exists\, 1 \le k \neq l \le n,\, a_{k}'(x) \ge \phi(n),\, a_{l}'(x) \ge \phi(n)
	~ \text{for i.m. } n \in \mathbb{N} \big\}.\]
	It follows from Eq. \eqref{0-1-[20]} that $a'_n(x)$ is unbounded for Lebesgue almost every \(x\). If $\phi(n)$ is bounded, such $x$ belongs to $E'(\phi)$. Thus, in this case, $\mathcal{L}(E'(\phi))=1$. Now, we assume that \( \lim_{n \to \infty} \phi(n) = \infty \).
	%	Without loss of generality, we may assume that \( \lim_{n \to \infty} \phi(n) = \infty \).
	%	Conversely, by Eq. \eqref{0-1-[20]}, $a'_n(x)$ is unbounded for Lebesgue almost every \(x\). When $\phi(n)$ is bounded, such $x$ belongs to $E'(\phi)$. Thus, in this case, $\mathcal{L}(E'(\phi))=1$. 
	
	We observe that \( E'(\phi) =\tilde{E}(\phi)\) where 
	\begin{equation}\label{eq:Ephi}
		\tilde{E}(\phi):=\big\{ x \in [0,1) : \exists\, 1 \le k < n,\, a_{k}'(x) \ge \phi(n),\, a_{n}'(x) \ge \phi(n)
		~ \text{for i.m.}\ n \in \mathbb{N} \big\}.
	\end{equation}
	Clearly, $\tilde{E}(\phi)\subseteq E'(\phi)$. To prove the reverse inclusion, for any \( x \in E'(\phi) \), there exists an increasing sequence \( (n_t)_{t \ge 1} \) such that for all \( t \ge 1 \),
	\[
	\exists\, 1 \le k_t < l_t \le n_t
	\quad \text{with} \quad
	a_{k_t}'(x) \ge \phi(n_t), \quad a_{l_t}'(x) \ge \phi(n_t).
	\]
	Since \( \lim_{n \to \infty} \phi(n) = \infty \), we see that $l_t$ tends to infinity as $t\to \infty$.
	It follows from the monotonicity of \( \phi \) that for all \( t\geq 1\),
	\[
	a_{k_t}'(x) \ge \phi(n_t) \ge \phi(l_t),
	\qquad
	a_{l_t}'(x) \ge \phi(n_t) \ge \phi(l_t).
	\]
	Therefore, $x$ belongs to $\tilde{E}(\phi)$. Henceforth, we adopt \eqref{eq:Ephi} as the definition of $E'(\phi)$.
	
	We first establish an auxiliary lemma for estimating $\mathcal{L}(E'(\phi))$.
	\begin{lem}\label{lem:mixing}
		Let $E_{k,n}:=\{x\in[0,1)\colon a'_k(x)\geq \phi(n)\}$ where  $1\leq k\leq n$. Let $F=\cup_{\mathbf{a}\in\mathcal{A}}I_{k-1}(\mathbf{a})$ where $\mathcal{A}\subset \mathbb{N}^{k-1}$.  Then 
		\[\mathcal{L}(F\cap E_{k,n})\asymp\frac{\mathcal{L}(F)}{\phi(n)\log\phi(n)}.\]
	\end{lem}
	\begin{proof}
		By Lemma \ref{lemma2.11} and Proposition \ref{p2.8}, we have
		\begin{align*}
			\mathcal{L}(F\cap E_{k,n}) & = \mathcal{L}(\cup_{\mathbf{a}\in\mathcal{A}}I_{k-1}(\mathbf{a})\cap E_{k,n})\\
			%\mathcal{L}\{x\in F\colon a'_k(x)\geq \phi(n)\} \\
			& = \sum_{\mathbf{a}\in\mathcal{A}}\mathcal{L}\{x\in I_{k-1}(\mathbf{a})\colon a'_k(x)\geq \phi(n)\}\\
			& = \sum_{\mathbf{a}\in \mathcal{A}}\sum_{a'_k\geq \phi(n)}\mathcal{L}\big(I_k(\mathbf{a}, a_{k})\big) = \sum_{\mathbf{a}\in \mathcal{A}}\sum_{p\in\mathcal{P}\atop p\geq \phi(n)}\mathcal{L}\big(I_k(\mathbf{a}, p)\big)\\
			& \asymp \sum_{\mathbf{a}\in \mathcal{A}}\sum_{p\in \mathcal{P},\atop p\geq \phi(n)}\frac{\mathcal{L}(I_{k-1}(\mathbf{a}))}{p^2} = \sum_{\mathbf{a}\in \mathcal{A}}\mathcal{L}(I_{k-1}(\mathbf{a}))\sum_{p\in \mathcal{P},\atop p\geq \phi(n)}\frac{1}{p^2}\\
			& \asymp  \frac{1}{\phi(n)\log\phi(n)}\sum_{\mathbf{a}\in \mathcal{A}}\mathcal{L}\big(I_{k-1}(\mathbf{a})\big) = \frac{\mathcal{L}(F)}{\phi(n)\log\phi(n)}.\qedhere
		\end{align*}
	\end{proof}

	\subsection{The convergence part}
	Note that \[E'(\phi)=\limsup_{n\to\infty} E'_n\] where
	$E'_n=\{x\in[0,1)\colon~\exists\ 1 \leq k< n,\ a_{k}'(x)\geq\phi(n),\ a_{n}'(x)\geq\phi(n)\}$. Write \[E_{k,n}:=\{x\in[0,1)\colon a'_k(x)\geq \phi(n)\}.\] Then $E'_n=\cup_{k=1}^{n}E_{k,n}\cap E_{n,n}$. 
	% We have 
	% \begin{align}
		% 	\mathcal{L}(E_{k,n}\cap E_{n,n}) & \leq \sum_{\mathbf{a}\in \mathbb{N}^{n-1}\atop a'_k\geq \phi(n)}\sum_{a'_n\geq \phi(n)}\mathcal{L}\big(I(\mathbf{a}, a_{n})\big) \lesssim \sum_{\mathbf{a}\in \mathbb{N}^{n-1}\atop a'_k\geq \phi(n)}\mathcal{L}\big(I(\mathbf{a})\big)\sum_{p\in \mathcal{P},\atop p\geq \phi(n)}\frac{1}{p^2}\nonumber\\
		% 	& \lesssim  \frac{1}{\phi(n)\log\phi(n)}\sum_{\mathbf{a}\in \mathbb{N}^{n-1}\atop a'_k\geq \phi(n)}\mathcal{L}\big(I(\mathbf{a})\big) = \frac{1}{\phi(n)\log\phi(n)}\mathcal{L}(E_{k,n}) \label{eq:3-1}
		% \end{align}
	% where two `$\lesssim$' follows from Lemma \ref{lemma2.11} and Proposition \ref{p2.8}. Similarly, 
	% \begin{equation}
		% 	\mathcal{L}(E_{k,n})\lesssim \frac{1}{\phi(n)\log\phi(n)}\sum_{\mathbf{a}\in\mathbb{N}^{k-1}}\mathcal{L}(I(\mathbf{a}))=\frac{1}{\phi(n)\log\phi(n)}.\label{eq:3-2}
		% \end{equation}	
	By Lemma \ref{lem:mixing}, \[\mathcal{L}(E'_n)\leq \sum_{k=1}^{n}\mathcal{L}(E_{k,n}\cap E_{n,n})\lesssim \sum_{k=1}^{n}\frac{1}{\phi^2(n)\log^2\phi(n)} = \frac{n}{\phi^2(n)\log^2\phi(n)}.\]
	Therefore, if $\sum^{\infty}_{n=1} \frac{n}{\phi^2(n)\log^2\phi(n)}<\infty$, then $\sum^{\infty}_{n=1}\mathcal{L}(E'_n)<\infty$. By Borel-Cantelli Lemma, we have $\mathcal{L}(E'(\phi))=0$.
	
	%	We have
	%	\begin{align*}
		%		\mathcal{L}(E'_n)%&\leq \sum_{1\leq k<n}\sum_{(a_1,\dots,a_{n-1})\in \mathbb{N},\atop a'_k\geq \phi(n),a'_n\geq \phi(n)}\mathcal{L}\big(I(a_{1}, a_{2}, \ldots, a_{n})\big)\\
		%		&\leq\sum_{1\leq k<n}\sum_{\mathbf{a}\in \mathbb{N}^{n-1}\atop a'_k\geq \phi(n)}\sum_{a'_n\geq \phi(n)}\mathcal{L}\big(I(\mathbf{a}, a_{n})\big)\\
		%		&\lesssim \sum_{1\leq k<n}\sum_{\mathbf{a}\in \mathbb{N}^{n-1}\atop a'_k\geq \phi(n)}\mathcal{L}\big(I(\mathbf{a})\big)\sum_{p\in \mathcal{P},\atop p\geq \phi(n)}\frac{1}{p^2}\qquad\quad \text{(by Lemma \ref{lemma2.11})}\\
		%		&\lesssim\sum_{1\leq k<n}\sum_{\mathbf{a}\in \mathbb{N}^{n-1}\atop a'_k\geq \phi(n)}\mathcal{L}\big(I(\mathbf{a})\big)\frac{1}{\phi(n)\log\phi(n)}.\quad \text{(by Proposition \ref{p2.8})}
		%	\end{align*}
	%	Applying Lemma \ref{lemma2.11} and Proposition \ref{p2.8} again, we have
	%	\begin{align*}
		%		\mathcal{L}(E'_n)&\lesssim\sum_{1\leq k<n}\sum_{(a_1,\dots,a_{k-1})\in \mathbb{N}}\sum_{a_k=p\in \mathcal{P},\atop p\geq \phi(n)}\mathcal{L}\big(I(a_{1}, a_{2}, \ldots, a_{k})\big)\frac{1}{\phi(n)\log\phi(n)}\\
		%		&\lesssim\sum_{1\leq k<n} \frac{1}{\phi(n)^2\log^2\phi(n)}\leq \frac{n}{\phi(n)^2\log^2\phi(n)}.
		%	\end{align*}
	%	Finally, Borel--Cantelli Lemma implies
	%	$$\mathcal{L}(E'(\phi))=0\quad\text{when} \quad \sum^{\infty}_{n=1} \frac{n}{\phi^2(n)\log^2\phi(n)}<\infty.$$
	
	\subsection{The divergence part}
	
	% Now suppose $\sum^{\infty}_{n=1}\frac{n}{\phi^2(n)\log^2\phi(n)}=\infty$. Without loss of generality, assume that $\phi(n)\geq \frac{n}{\log \phi(n)}$ for all $n\geq 1$. Otherwise, we define
	% $$\hat{\phi}(n)=\max\{\phi(n),\ \tfrac{n}{\log \phi(n)}\}.$$
	% Then $E'(\hat{\phi})\subset E'(\phi)$, and the divergence still holds:
	% \begin{equation*}
	% 	\sum^{\infty}_{n=1}\frac{n}{\phi^2(n)\log^2\phi(n)}=\infty \implies \sum^{\infty}_{n=1}\frac{n}{\hat{\phi}^2(n)\log^2\hat{\phi}(n)}=\infty.
	% \end{equation*}
	% In fact, let $\mathcal{N}=\{n\in\mathbb{N}\colon \phi(n)<\frac{n}{\log \phi(n)}\}$. If $\#\mathcal{N}<\infty$, then $\hat\phi$ is ultimately $\phi$ and the series diverges. In the case $\#\mathcal{N}=\infty$, choose $(n_k)_{k\geq 1}\subset \mathcal{N}$ with $n_{k+1}>\frac{1}{c}\cdot n_{k}$ for some $c\in (0,1)$. Then $\phi(n_k)<\frac{n_k}{\log \phi(n_k)}$ and $\hat\phi(n_k)=\frac{n_k}{\log \phi(n_k)}$ for all $k\geq 1$. Consequently, 
	% \begin{align*}
	% 	\sum^{\infty}_{n=1}\frac{n}{\hat{\phi}^2(n)\log^2\hat{\phi}(n)} & \geq \sum^{\infty}_{k=1}\sum_{cn_{k+1}\leq n\leq n_{k+1}}\frac{n}{\hat{\phi}^2(n)\log^2\hat{\phi}(n)} \\
	% 	& \geq  \sum^{\infty}_{k=1}\sum_{cn_{k+1}\leq n\leq n_{k+1}}\frac{cn_{k+1}}{\hat{\phi}^2(n_{k+1})\log^2\hat{\phi}(n_{k+1})}\\
	% 	& \geq \sum^{\infty}_{k=1}\frac{c(1-c)n^2_{k+1}}{n^2_{k+1}} \geq \sum^{\infty}_{k=1}c(1-c) = \infty.
	% \end{align*}
	
	Now suppose $\sum^{\infty}_{n=1}\frac{n}{\phi^2(n)\log^2\phi(n)}=\infty$. Without loss of generality, assume that $\phi(n)\log\phi(n)\geq n$ for all $n\geq 1$. Otherwise, we define
		$$\hat{\phi}(n)=\max\{\phi(n),\ \tfrac{2n}{\log n}\}.$$
		Then $E'(\hat{\phi})\subset E'(\phi)$, and the divergence still holds:
		\begin{equation*}
			\sum^{\infty}_{n=1}\frac{n}{\phi^2(n)\log^2\phi(n)}=\infty \implies \sum^{\infty}_{n=1}\frac{n}{\hat{\phi}^2(n)\log^2\hat{\phi}(n)}=\infty.
		\end{equation*}
		It is clearly that $\hat\phi(n)$ is increasing. Moreover, for all $n\geq 2$,
		\[\log\bigl(\frac{2n}{\log n}\bigr)=\log n+\log2-\log\log n\geq \frac{1}{2} \log n.\] 
		Therefore, for all $n\geq 2$, we have
		\[\hat{\phi}(n)\log\hat{\phi}(n)\geq \frac{2n}{\log n}\log(\frac{2n}{\log n})\geq n.\]
		In fact, let $\mathcal{N}=\{n\in\mathbb{N}\colon \phi(n)<\frac{2n}{\log n}\}$. If $\#\mathcal{N}<\infty$, then $\hat\phi$ is ultimately $\phi$ and the series diverges. In the case $\#\mathcal{N}=\infty$, choose $(n_k)_{k\geq 1}\subset \mathcal{N}$ with $n_{k+1}>\frac{1}{c}\cdot n_{k}$ for some $c\in (0,1)$. Then $\phi(n_k)<\frac{2n_k}{\log n_k}$ and $\hat\phi(n_k)=\frac{2n_k}{\log n_k}$ for all $k\geq 1$. Consequently, 
		\begin{align*}
			\sum^{\infty}_{n=1}\frac{n}{\hat{\phi}^2(n)\log^2\hat{\phi}(n)} & \geq \sum^{\infty}_{k=1}\sum_{cn_{k+1}\leq n\leq n_{k+1}}\frac{n}{\hat{\phi}^2(n)\log^2\hat{\phi}(n)} \\
			& \geq  \sum^{\infty}_{k=1}\sum_{cn_{k+1}\leq n\leq n_{k+1}}\frac{cn_{k+1}}{\hat{\phi}^2(n_{k+1})\log^2\hat{\phi}(n_{k+1})}\\
			& \geq \sum^{\infty}_{k=1}\frac{c(1-c)n^2_{k+1}}{n^2_{k+1}} \geq \sum^{\infty}_{k=1}c(1-c) = \infty.
	\end{align*}

	Fix $t\geq 1$ and $\mathbf{a}=(a_1,\dots,a_t)\in\mathbb{N}^{t}$, let $I_t:=I_t(\mathbf{a})$. We estimate $\mathcal{L}(E'(\phi)\cap I_t)$.
	Recall that $E'(\phi)\cap I_t=\limsup_{n\rightarrow \infty}E'_n\cap I_t$. Choose $n_0\in \mathbb{N}$ such that $n\geq 2t+2$ and $\phi(n)\geq \max\{a_1,\dots,a_t\}$ for all $n\geq n_0$. Write $F_k=\{x\in I_t\colon a'_i(x)<\phi(n)\  \text{for all} \ t<i<k\}$. Then \[E'_n\cap I_t = \bigcup_{k=t+1}^{n-1}F_k\cap E_{k,n}\cap E_{n,n}.\] Note that the above union is disjoint. By Lemma \ref{lem:mixing}, 
	% \[\mathcal{L}(F_k\cap E_{k,n}\cap E_{n,n}) \gtrsim \frac{1}{\phi(n)^2\log^2\phi(n)}\mathcal{L}(F_k)\]
	\[\mathcal{L}(E'_n\cap I_t) = \sum^{n-1}_{k=t+1}\mathcal{L}(F_k\cap E_{k,n} \cap E_{n,n}) \gtrsim \sum^{n-1}_{k=t+1}\frac{1}{\phi(n)^2\log^2\phi(n)}\mathcal{L}(F_k).\]
	Note that \[F_k=I_t\cap \bigcap_{i=t+1}^{k-1}E^{c}_{i,n}.\] Using Lemma \ref{lem:mixing} several times, there exist $c_{n_0}> 0$ for all $n\geq n_0$, 
	\[\mathcal{L}(F_k)\gtrsim \mathcal{L}(I_t)\Bigl(1-\frac{c_{n_0}}{\phi(n)\log\phi(n)}\Bigr)^{k-t-1}\geq \mathcal{L}(I_t)\Bigl(1-\frac{c_{n_0}}{n}\Bigr)^{k-t-1}\] where the last inequality holds since $\phi(n)\geq \frac{n}{\log \phi(n)}$.
	We conclude that 
	\begin{align}
		\mathcal{L}(E'_n\cap I_t) & \gtrsim \sum^{n-1}_{k=t+1}\frac{\mathcal{L}(F_k)}{\phi(n)^2\log^2\phi(n)}\nonumber\\ 
		& \gtrsim \frac{\mathcal{L}(I_t)}{\phi(n)^2\log^2\phi(n)}\sum^{n-1}_{k=t+1}\left(1-\frac{c_{n_0}}{n}\right)^{k-t-1}\nonumber\\
		& \gtrsim \frac{n}{\phi(n)^2\log^2\phi(n)}\mathcal{L}(I_t).\label{eq:step-1}
	\end{align}
	
	We estimate the measure of $\mathcal{L}((E'_n\cap I_t)\cap(E'_m\cap I_t))$ for $n>m\geq n_0$.
	
	Observe that
	\begin{align*}
		&(E'_n\cap I_t)\cap (E'_m\cap I_t)\\
		= & \{x\in I_t: \exists\ l<m,\ k<n, a'_l(x)\geq\phi(m), a'_m(x)\geq\phi(m), a'_k(x)\geq \phi(n), a'_n(x)\geq \phi(n)\}\\
		= & \bigcup_{k=t+1}^{n-1}\bigcup_{l=t+1}^{m-1}\bigl(I_{t}\cap E_{l,m}\cap E_{m,m}\cap E_{k,n}\cap E_{n,n}\bigr).
	\end{align*}
	When $k=l$, by the monotonicity of $\phi$, we have $E_{k,n}\subset E_{l,m}$ and by Lemma \ref{lem:mixing}, 
	\begin{align*}
		\mathcal{L}\bigl(I_{t}\cap E_{l,m}\cap E_{m,m}\cap E_{k,n}\cap E_{n,n}\bigr) & = \mathcal{L}\bigl(I_{t}\cap E_{m,m}\cap E_{k,n}\cap E_{n,n}\bigr)\\
		& \lesssim \frac{1}{\phi(m)\log\phi(m)}\cdot \frac{1}{\phi^2(n)\log^2\phi(n)}\cdot\mathcal{L}(I_t).
	\end{align*}
	When $k=m$, we also have $E_{k,n}\subset E_{m,m}$ and by Lemma \ref{lem:mixing}, 
	\begin{align*}
		\mathcal{L}\bigl(I_{t}\cap E_{l,m}\cap E_{m,m}\cap E_{k,n}\cap E_{n,n}\bigr) & = \mathcal{L}\bigl(I_{t}\cap E_{l,m}\cap E_{k,n}\cap E_{n,n}\bigr)\\
		& \lesssim \frac{1}{\phi(m)\log\phi(m)}\cdot \frac{1}{\phi^2(n)\log^2\phi(n)}\cdot\mathcal{L}(I_t).
	\end{align*}
	When $k\neq l$ or $m$,  by Lemma \ref{lem:mixing},  we have 
	\begin{align*}
		\mathcal{L}\bigl(I_{t}\cap E_{l,m}\cap E_{m,m}\cap E_{k,n}\cap E_{n,n}\bigr) & \lesssim \frac{1}{\phi^2(m)\log^2\phi(m)}\cdot \frac{1}{\phi^2(n)\log^2\phi(n)}\cdot\mathcal{L}(I_t).
	\end{align*}
	Therefore, 
	\begin{align}
		\mathcal{L}\bigl((E'_n\cap I_t)\cap (E'_m\cap I_t)\bigr) & \leq \sum_{k=t+1}^{n-1}\sum_{l=t+1}^{m-1}\mathcal{L}\bigl(I_{t}\cap E_{l,m}\cap E_{m,m}\cap E_{k,n}\cap E_{n,n}\bigr)\notag\\
		& \lesssim  \frac{2(m-t-1)}{\phi(m)\log\phi(m)}\cdot\frac{1}{\phi^2(n)\log^2\phi(n)}\cdot\mathcal{L}(I_t)\notag\\
		& \qquad + \frac{(m-t-1)(n-t-3)}{\phi^2(m)\log^2\phi(m)}\cdot \frac{1}{\phi^2(n)\log^2\phi(n)}\cdot\mathcal{L}(I_t)\notag\\
		& \lesssim \frac{1}{\phi^2(n)\log^2\phi(n)}\cdot\mathcal{L}(I_t) + \mathcal{L}(E'_m\cap I_t)\cdot\mathcal{L}(E'_n\cap I_t)\cdot\frac{1}{\mathcal{L}(I_t)}\label{eq:step-2}
	\end{align}
	where in the last step we need $\phi(m)\geq \frac{m}{\log\phi(m)}$ and Eq. \eqref{eq:step-1}.
	
	Using Eq. \eqref{eq:step-2} and Eq. \eqref{eq:step-1}, we have 
	\begin{align*}
		& \sum_{n_0\leq m< n\leq N}\mathcal{L}((E'_n\cap I_t)\cap(E'_m\cap I_t))\\
		\lesssim & \sum_{n=n_0}^{N}\sum_{m=n_0}^{n-1}\Bigl(\frac{1}{\phi^2(n)\log^2\phi(n)}\cdot\mathcal{L}(I_t) + \mathcal{L}(E'_m\cap I_t)\cdot\mathcal{L}(E'_n\cap I_t)\cdot\frac{1}{\mathcal{L}(I_t)}\Bigr)\\
		= & \sum_{n=n_0}^{N}\frac{n-n_0}{\phi^2(n)\log^2\phi(n)}\cdot\mathcal{L}(I_t) + \sum_{n=n_0}^{N}\sum_{m=n_0}^{n-1}\mathcal{L}(E'_m\cap I_t)\cdot\mathcal{L}(E'_n\cap I_t)\cdot\frac{1}{\mathcal{L}(I_t)}\\
		\lesssim & \sum_{n=n_0}^{N}\mathcal{L}(E'_n\cap I_t) + \Bigl(\sum_{n=n_0}^{N}\mathcal{L}(E'_n\cap I_t)\Bigr)^2\cdot\frac{1}{\mathcal{L}(I_t)}.
	\end{align*}
	Consequently, since $\sum^{\infty}_{n=1}\mathcal{L}(E'_n\cap I_t)=\infty$, we have 
	\begin{align*}
		\frac{\Bigl(\sum_{n=n_0}^{N}\mathcal{L}(E'_n\cap I_t)\Bigr)^2}{\sum\limits_{n_0\leq m, n\leq N}\mathcal{L}(E'_m\cap E'_n\cap I_t)} & \gtrsim \frac{\Bigl(\sum_{n=n_0}^{N}\mathcal{L}(E'_n\cap I_t)\Bigr)^2}{\sum_{n=n_0}^{N}\mathcal{L}(E'_n\cap I_t) + \left(\sum_{n=n_0}^{N}\mathcal{L}(E'_n\cap I_t)\right)^2\cdot\frac{1}{\mathcal{L}(I_t)}}\to \mathcal{L}(I_t)
	\end{align*}
	as $N\to\infty$. Applying Chung-Erd\"os inequality (Lemma~\ref{L2.15}), we obtain
	\[
	\mathcal{L}(E'(\phi)\cap I_t)
	=\mathcal{L}\bigr(\limsup_{n\to\infty} E'_n\cap I_t\bigr)
	\gtrsim \mathcal{L}(I_t).
	\]
	Finally, by Lemma~\ref{L2.16}, we conclude that $\mathcal{L}(E'(\phi))=1$.

	\section{Proof of Theorem \ref{t2}: the case $1<B<\infty$}

	In this section, we calculate the Hausdorff dimension of $E'(\phi)$ for the case $1<B<\infty$. Recall that $\log B=\liminf\limits_{n\to \infty}\frac{\log\phi(n)}{n}$ and $T$ is the Gauss map. Consider the potential
	\begin{equation*}%\label{eq:phi}
		f(x) = -(3s-1)\log B - s\log|T'(x)|.
	\end{equation*}
	For $A\subseteq\mathbb{N}$, define
	\[
	s(B,A) := \inf\{ s \ge 0 : P_A(T, f) \le 0 \}.
	\]
	Write $s(B,M):=s(B,\{1,\dots,M\})$ where $M\geq 1$ and $s(B):=s(B,\mathbb{N})$.  
	Since $E'(\phi) \subset E(\phi)$, we have \[\dim_{H}E'(\phi)\leq \dim_{H}E(\phi)=s(B)\] where the last equality is given in \cite[Theorem 1.6]{TTW}.  
	
	\subsection{Construction of a Cantor Set} 	
	To obtain the lower bound, we construct a Cantor-type subset together with a suitable mass distribution.
	
	\subsubsection{Setting of parameters}\label{sec:4.1}
	For all $n\geq 1$, let $s_n(B,A)$ be the unique solution of
	\begin{equation*}%\label{eq:sn}
		\sum_{(a_1,\dots,a_n)\in A^n} \frac{1}{B^{(3s-1)n} q_n^{2s}(a_1,\dots,a_n)} = 1.
	\end{equation*}
	Write $s_n(B,M):=s_n(B,\{1,\dots,M\})$ for $M\geq 1$ and $s_n(B):=s_n(B,\mathbb{N})$. It is known that the functions $s_n(\cdot)$, $s_n(\cdot,\cdot)$, $s(\cdot)$ and $s(\cdot,\cdot)$ satisfy the following properties.
	
	% Recall that $\log B=\liminf\limits_{n\to \infty}\frac{\log\phi(n)}{n}$ and $T$ is the Gauss mapping. Consider the potential
	% \begin{equation*}%\label{eq:phi}
		% 	f(x) = -(3s-1)\log B - s\log|T'(x)|.
		% \end{equation*}
	% For $A\subseteq\mathbb{N}$, define
	% \[
	% s(B,A) := \inf\{ s \ge 0 : P_A(T, f) \le 0 \},
	% \]
	% and let $s_n(B,A)$ be the unique solution of
	% \begin{equation*}%\label{eq:sn}
		% 	\sum_{(a_1,\dots,a_n)\in A^n} \frac{1}{B^{(3s-1)n} q_n^{2s}(a_1,\dots,a_n)} = 1.
		% \end{equation*}
	% For any fixed $M\in\mathbb{N}$ and $A=\{1,2,\dots,M\}$, write $s_n(B,M):=s_n(B,A)$ and $s(B,M):=s(B,A)$.
	% Similarly, we write $s_n(B):=s_n(B,\mathbb{N})$ and $s(B):=s(B,\mathbb{N})$.
	% We have the following properties:
	\begin{lem}[{\cite[Proposition 2.11]{TTW}}]\label{pro4.1}
		The following assertions hold:
		\begin{enumerate}
			\item For all $n\geq 1$ and $B>1$, $s_n(B,\mathbb{N})\ge \tfrac{1}{2}$;
			\item $\lim\limits_{n\to\infty} s_n(B,M) = s(B,M)$, $\lim\limits_{n\to\infty} s_n(B) = s(B)$, and $\lim\limits_{M\to\infty} s(B,M) = s(B)$;
			\item As a function of $B\in(1,\infty)$, $s(B)$ is strictly decreasing and continuous;
			\item $\lim\limits_{B\to1} s(B) = 1$ and $\lim\limits_{B\to\infty} s(B) = \tfrac{1}{2}$.
		\end{enumerate}
	\end{lem}
	
	By the definition of $B$, let $\{n_k\}_{k\ge 1}$ be a sequence of integers such that
	\[
	\log B = \lim_{k\to\infty} \frac{\log\phi(n_k + 2)}{n_k + 2}.
	\]
	Let $s\in \bigl(\tfrac{1}{2},s(B)\bigr)$. By Lemma \ref{pro4.1}, we can choose $\tilde{B} > B$ such that $s < s(\tilde{B}) < s(B)$ and $\phi(n_{k}+2)\leq \tilde{B}^{n_k}$ for all $k$, and for all large enough $M$ and $N$, one has 
	\begin{equation*}%\label{eq4.2}
		s_N(\tilde{B},M) > s.
	\end{equation*}
	Let $\delta\in (0, s-\tfrac{1}{2})$. We may assume that $N\geq \max\{e^{20},\tfrac{2}{\delta}+1\}$.
	
	Let $\{\ell_k\}_{k\ge 1}$ be the sequence of integers satisfying for all $k\geq 0$, 
	\[
	n_{k+1} - n_k = 2 + \ell_{k+1} N + N + i_{k+1}
	\]
	where $0 \le i_{k+1} < N$ and $n_0:=-2$. Thus $N\sum_{j=1}^{k}\ell_j\leq n_k\leq N\sum_{j=1}^{k}(\ell_j +3)$ for all $k\geq 1$. Moreover, let
	\[
	m_k := n_k-N-i_k= n_{k-1} + 2 + \ell_k N.
	\]
	By choosing subsequences of $\{n_k\}_{k\geq 1}$, we may assume that $\{\ell_k\}_{k\geq 1}$ increases rapidly, for example $\ell_k\geq e^{k^2}$, so that
	\begin{equation}\label{eq4.3}
		\frac{N \cdot \ell_k}{n_k} \ge 1 - \delta, \quad
		\frac{\log n_k}{n_k} < \delta, \quad \text{and} \quad
		\ell_k \geq 24N \quad \text{for all } k \ge 1,
	\end{equation}

	% \color{red}
	%  Let $s$ and $\delta$ be positive real numbers such that
	% \[
	% \frac{1}{2} < s < s(B) \quad \text{and} \quad s - \delta > \frac{1}{2}.
	% \]
	% Hence, it is enough to show that $\dim_H E'(\psi) \ge s$ for all  $s < s(B)$.

	% Choose $\tilde{B} > B$ such that $s < s(\tilde{B}) < s(B)$, and fix $M$ and $N$  large enough so that
	% \begin{equation}\label{eq4.2}
		% 	s_N(\tilde{B},M) > s, \quad \text{and} \quad N > \max\{ e^{20}, \tfrac{2}{\delta} + 1 \}.
		% \end{equation}
	% By the definition of $B$, let $\{n_k\}_{k\ge 1}$ be the sequence of integers such that
	% \[
	% \log B = \lim_{k\to\infty} \frac{\log\phi(n_k + 2)}{n_k + 2}, \quad \phi(n_k + 2) \leq \tilde{B}^{\,n_k}\ \text{for all}\ k \ge 1.
	% \]
	
	% Let $\{\ell_k\}_{k\ge 1}$ be a rapidly increasing sequence of integers and define
	% \[
	% n_1 = \ell_1 N + N + i_1, \qquad
	% n_{k+1} - n_k = 2 + \ell_{k+1} N + N + i_{k+1} \quad \text{for all } k \ge 1,
	% \]
	% where $0 \le i_k < N$. Moreover, set
	% \[
	% m_k := n_k-N-i_k= n_{k-1} + 2 + \ell_k N, \quad \text{with } n_0 = -2.
	% \]
	
	% We assume that $\{n_k\}$ is sufficiently sparse so that
	% \begin{equation}\label{eq4.3}
		% 	\frac{N \cdot \ell_k}{n_k} \ge 1 - \delta, \quad
		% 	\frac{\log n_k}{n_k} < \delta, \quad \text{and} \quad
		% 	\ell_k \geq 24N \quad \text{for all } k \ge 1,
		% \end{equation}

	\subsubsection{Constructing the Cantor set}
	Let $E_{M}^{N}(\tilde{B})$ be the set of irrational numbers $x \in [0,1)$ satisfying the following conditions:
	\begin{enumerate}
		\item $a'_{n_j+1}(x),\ a'_{n_j+2}(x) \in [\tilde{B}^{n_j},\, 2\tilde{B}^{n_j}]$ for all $j\geq 1$;
		\item $a_{m_j+1}(x)=\dots=a_{n_j}(x)=2$ for all $j\geq 1$;
		\item  $1 \leq a_n(x) \leq M$ for all remaining indices $n\geq 1$.
	\end{enumerate}
	That is, for any $x\in E_{M}^{N}(\tilde{B})$, the partial quotients of $x$ satisfy
	\[\bigl(\,\underbrace{a_1(x), \dots, a_{m_1}(x)}_{\in \{1,\dots,M\}^{\ell_1 N}},\ \underbrace{a_{m_1+1}(x),\dots, a_{n_1}(x)}_{N+i_1\text{ repetitions of }2},\ \underbrace{a_{n_1+1}(x),a_{n_1+2}(x)}_{\in [\tilde{B}^{n_1},\, 2\tilde{B}^{n_1}]^2\cap\mathcal{P}^2},\ \underbrace{a_{n_1+3}(x), \dots, a_{m_2}(x)}_{\in \{1,\dots,M\}^{\ell_2 N}}, \dots \,\bigr)\]
	
	Since $\tilde{B}^{n_k}\geq \phi(n_k+2)$ and $\phi$ is non-decreasing, we have $E_{M}^{N}(\tilde{B}) \subset E'(\phi)$.
	In what follows, we describe this set in terms of cylinder sets.
	For all $n \in \mathbb{N}$, define
	\[
	D_n := \{ \mathbf{a}\in A^n \colon I_n(\mathbf{a}) \subset E_{M}^{N}(\tilde{B}) \}\quad \text{and}\quad 
	J_n(\mathbf{a}) := \bigcup_{(\mathbf{a}, a_{n+1}) \in D_{n+1}} I_{n+1}(\mathbf{a},a_{n+1}).
	\]
	Then 
	\[
	E_{M}^{N}(\tilde{B}) = \bigcap_{n \geq 1} \bigcup_{\mathbf{a}\in D_n} J_n(\mathbf{a}).
	\]
	\subsubsection{The lengths of $J_n(\mathbf{a})$}
	
	% In this section, we estimate the lengths of the cylinders
	% $\{J_n(\mathbf{a})\}$, where $(\mathbf{a}) \in D_n$.
	Let $\mathbf{a}$ be the continued fraction expansion of some $x\in E_{M}^{N}(\tilde{B})$. Assume that $n_{j-1} + 2 \leq n \leq  n_j + 1$ for some $j \ge 1$.
	
	When $n_{j-1} + 2 \le n < m_j$, we have
	\begin{equation*}
		|J_n(\mathbf{a})| = \Bigl|\bigcup_{a_{n+1}=1}^{M}I_{n+1}(\mathbf{a},a_{n+1})\Bigr|=\sum_{a_{n+1}=1}^{M}|I_{n+1}(\mathbf{a},a_{n+1})|.
	\end{equation*}
	Since $\tfrac{1}{2}q_{n+1}^{-2}(\mathbf{a},a_{n+1})\leq |I_{n+1}(\mathbf{a},a_{n+1})|\leq q_{n+1}^{-2}(\mathbf{a},a_{n+1})$ by Lemma \ref{L2.9} and $a_{n+1}q_n(\mathbf{a})< q_{n+1}(\mathbf{a},a_{n+1})< (a_{n+1}+1)q_n(\mathbf{a})$, 
	\begin{equation}\label{eqjn1}
		\frac{1}{8 q_n^2(\mathbf{a})}\leq \frac{1}{q_n^2(\mathbf{a})}\sum_{a_{n+1}=1}^{M}\frac{1}{(a_{n+1}+1)^2}\leq |J_n(\mathbf{a})| \leq \frac{1}{q_n^2(\mathbf{a})}\sum_{a_{n+1}=1}^{M}\frac{1}{a_{n+1}^2}\leq \frac{2}{q_n^2(\mathbf{a})}.
	\end{equation}
	
	When $m_j \le n < n_j$, similarly, we obtain
	\begin{equation}\label{eqjn2}
		|J_n(\mathbf{a})|=|I_{n+1}(\mathbf{a},2)| \in \left[\frac{1}{18q_n^2(\mathbf{a})},\ \frac{1}{4q_n^2(\mathbf{a})}\right].
	\end{equation}
	
	Finally, when $n = n_j$ or $n=n_j+1$, we obtain
	\begin{align*}
		|J_n(\mathbf{a})|
		&= \Bigl|
		\bigcup_{a_{n+1}\in [\tilde{B}^{n_j}, 2 \tilde{B}^{n_j}]\cap \mathcal{P}}
		I_{n+1}(\mathbf{a}, a_{n+1})
		\Bigr| \\
		&= \biggl|
		\frac{p_{\text{min}}p_n(\mathbf{a})+p_{n-1}(\mathbf{a})}{p_{\text{min}}q_n(\mathbf{a})+q_{n-1}(\mathbf{a})}-\frac{p_{\text{max}}p_n(\mathbf{a})+p_{n-1}(\mathbf{a})}{p_{\text{max}}q_n(\mathbf{a})+q_{n-1}(\mathbf{a})}
		\biggr| \\
		& = \frac{(p_{\max}-p_{\min})}{(p_{\min}q_n(\mathbf{a})+q_{n-1}(\mathbf{a}))(p_{\max}q_n(\mathbf{a})+q_{n-1}(\mathbf{a}))}\\
		& \in \left[\frac{p_{\max}-p_{\min}}{(p_{\min}+1)(p_{\max}+1)q_n^2(\mathbf{a})},\ \frac{p_{\max}-p_{\min}}{p_{\min}p_{\max}q_n^2(\mathbf{a})}\right]
		% & \asymp \frac{p_{\text{max}}-p_{\text{min}}}{p_{\text{max}}p_{\text{min}}}\frac{1}{q_n^2(\mathbf{a})} \asymp \frac{p_{\text{max}}-p_{\text{min}}}{\tilde{B}^{2n_j}q_n^2(\mathbf{a})}
	\end{align*}
	where $p_{\text{min}}$ and $p_{\text{max}}$ are the smallest and the largest primes in $[\tilde{B}^{n_j}, 2 \tilde{B}^{n_j}]$. By Proposition~\ref{p2.5}~(2), we have 
	\[\tilde{B}^{n_j}\leq p_{\min}< \tfrac{1}{0.999}\tilde{B}^{n_j}\quad \text{and}\quad 0.999\cdot 2\tilde{B}^{n_j}\leq p_{\max}< 2\tilde{B}^{n_j}.\]
	It follows that
	\begin{equation}
		\frac{1}{3\tilde{B}^{n_j} q_n^2(\mathbf{a})}\leq |J_n(\mathbf{a})| \leq \frac{1}{\tilde{B}^{n_j} q_n^2(\mathbf{a})}. \label{eqjn3}
	\end{equation}
	% In fact, from the above, $|J_n(\mathbf{a})|$ is quite close to $\tfrac{1}{2\tilde{B}^{n_j} q_n^2(\mathbf{a})}$.

	% The fundamental set corresponding to $(\mathbf{a})$ is given by
	% \[
	% J_n(\mathbf{a}) = \bigcup_{(\mathbf{a}, a_{n+1}) \in D_{n+1}} I_{n+1}(\mathbf{a},a_{n+1}).
	% \]
	% Consequently,
	% \[
	% E_{M}^{N}(\tilde{B}) = \bigcap_{n \geq 1} \; \bigcup_{(\mathbf{a})\in D_n} J_n(\mathbf{a}).
	% \]

	\subsection{ Mass distribution on $E_{M}^{N}(\tilde{B})$}	
	Let $x\in E_{M}^{N}(\tilde{B})$ with 
	\[
	x = (w_1^{(1)}, \dots, w_{\ell_1}^{(1)}, v^{(1)}, a_{n_1+1}, a_{n_1+2},
	w_1^{(2)}, \dots, w_{\ell_2}^{(2)}, v^{(2)}, a_{n_2+1}, a_{n_2+2}, \dots \,),
	\]
	where $w_i^{(k)} \in \{1, \dots, M\}^{N}$ for all $i, k \in \mathbb{N}$. 
	
	Let $\mu(J_0)=1$ and \[\mu(I_{N}(w))=\frac{1}{u}\cdot\frac{1}{\tilde{B}^{(3s-1)N}q^{2s}_N(w)}\] for all $w\in A^N$ where \[u=\sum_{w\in \{1,\dots,M\}^{N}}\frac{1}{\tilde{B}^{(3s-1)N}q^{2s}_N(w)}.\] Now we define $\mu$ recursively. For all $k\geq 1$,  
	\begin{itemize}
		\item let $\mu(J_{m_k}(a_1,\dots,a_{m_k})) := \mu(J_{n_{k-1}+2}(a_1,\dots,a_{n_{k-1}+2}))\prod_{i=1}^{\ell_k}\mu(I_{N}(w_i^{(k)}))$;
		\item for $m_k< n \leq n_k$, define 
		\[\mu(J_{n}(a_1,\dots,a_n))=\mu(J_{m_k}(a_1,\dots,a_{m_k}));\] 
		%if $(a_{m_k+1},\dots,a_n)=(2,\dots,2)$, and $\mu(J_{n}(a_1,\dots,a_n))=0$ otherwise;
		\item for $n=n_k+i$ with $i=1$ or $2$, \[\mu(J_{n}(a_1,\dots,a_n))=\frac{1}{(\# P_{n_k})^i}\mu(J_{m_k}(a_1,\dots,a_{m_{k}}));\]
		\item for $n_{k-1}+2 < n < m_{k}$, we have $J_n(x)=\cup J_{m_k}(x)$ \[\mu(J_n(a_1,\dots,a_n))=\sum_{J_{m_k}(a_1, \dots, a_{m_k}) \subset J_n}\mu(J_{m_k}(a_1,\dots,a_{m_{k}})).\]
	\end{itemize}
	
	\subsubsection{H\"older exponent of fundamental intervals}
	From the construction of $\mu$, evaluating $\mu(J_n)$ reduces to estimating the product $\prod_{i=1}^{\ell_k}\mu(I_{N}(w_i^{(k)}))$, which is bounded in the lemma below.
	\begin{lem}\label{lem:prod-mu-i}
		Let $N$, $M$ and $\delta$ be as in Subsection \ref{sec:4.1}. Let $\ell\geq 1$ and $w_i\in\{1,\dots,M\}^N$ for all $i=1,\dots,\ell$. Then 
		\[\prod_{i=1}^{\ell}\mu(I_{N}(w_i))\leq \frac{1}{\tilde{B}^{(3s-1)N\ell}}\cdot\frac{1}{q^{2s(1-\delta)}_{\ell N}(w_1,\dots,w_\ell)}.\]
	\end{lem}
	\begin{proof}
		Since $u\geq 1$, it follows from the definition of $\mu(I_N(w_i))$ that 
		\begin{align*}
			\prod_{i=1}^{\ell}\mu(I_{N}(w_i)) & \leq \frac{1}{\tilde{B}^{(3s-1)N\ell}}\cdot\frac{1}{q^{2s}_N(w_1)\cdots q^{2s}_N(w_{\ell})} \leq \frac{1}{\tilde{B}^{(3s-1)N\ell}}\cdot\frac{4^{\ell s}}{q^{2s}_{\ell N}(w_1, \dots, w_{\ell})}
		\end{align*}
		where in the last step we apply Lemma \ref{eq:2}(2). Note that  $q_{\ell N}\geq 2^{(\ell N-1)/2}$ by Lemma \ref{eq:2}~(1). Recall that we have assumed $N\geq \tfrac{2}{\delta}+1$. So, $q_{\ell N}\geq 2^{\ell/\delta}$. Consequently, 
		\begin{align*}
			\prod_{i=1}^{\ell}\mu(I_{N}(w_i)) & \leq \frac{1}{\tilde{B}^{(3s-1)N\ell}}\cdot\frac{1}{q^{2s(1-\delta)}_{\ell N}(w_1, \dots, w_{\ell})}.\qedhere
		\end{align*}
	\end{proof}
	\medskip
	
	We now estimate $\mu(J_n)$.
	
	(1) For $n=n_{k-1}+2$, we have 
	\begin{align}
		\mu(J_{n_{k-1}+2}(a_1,\dots,a_{n_{k-1}+2})) & = \mu\bigl(J_{m_{k-1}}(a_1,\dots,a_{m_{k-1}})\bigr)\frac{1}{(\# P_{n_{k-1}})^2}\notag\\
		& = \mu\bigl(J_{n_{k-2}+2}(a_1,\dots,a_{n_{k-2}+2})\bigr)\frac{1}{(\# P_{n_{k-1}})^2}\prod_{i=1}^{\ell_{k-1}}\mu(I_{N}(w_i^{(k-1)}))\notag\\
		% & = \mu(J_{m_{k-1}}(a_1,\dots,a_{m_{k-1}}))\frac{1}{(\# P_{n_{k-1}})^2}\prod_{i=1}^{\ell_{k}}\mu(I_{N}(w_i^{(k)}))\\
		& = \prod_{t=1}^{k-1}\frac{1}{(\# P_{n_{t}})^2}\prod_{i=1}^{\ell_{t}}\mu(I_{N}(w_i^{(t)}))\notag\\
		(\text{by Proposition~\ref{p2.5}})\quad & = \prod_{t=1}^{k-1}\biggl(\frac{\log \tilde{B}^{n_t}}{c_{n_t}(\tilde{B})\tilde{B}^{n_t}}\biggr)^2\prod_{i=1}^{\ell_{t}}\mu(I_{N}(w_i^{(t)}))\notag\\
		& \leq \prod_{t=1}^{k-1}\frac{1}{\tilde{B}^{2n_t(1-\delta)}}\prod_{i=1}^{\ell_{t}}\mu(I_{N}(w_i^{(t)}))\notag\\
		(\text{by Lemma }\ref{lem:prod-mu-i})\quad & \leq \prod_{t=1}^{k-1}\frac{1}{\tilde{B}^{2n_t(1-\delta)}}\frac{1}{\tilde{B}^{(3s-1)N\ell_t}}\cdot\frac{1}{q^{2s(1-\delta)}_{\ell_t N}(w_1^{(t)}, \dots, w_{\ell_t}^{(t)})}\notag\\
		% (\text{Def of }\mu)\quad & = \prod_{t=1}^{k-1}\frac{1}{\tilde{B}^{2n_t(1-\delta)}} \frac{1}{u^{\ell_t}}\frac{1}{\tilde{B}^{(3s-1)N\ell_t}}\cdot\frac{1}{q^{2s}_N(w_1^{(t)})\cdots q^{2s}_N(w_{\ell_t}^{(t)})}\notag\\
		% (u\geq 1 \text{ and by Lemma }\ref{eq:2})\quad & \leq \prod_{t=1}^{k-1}\frac{1}{\tilde{B}^{2n_t(1-\delta)}} \frac{1}{\tilde{B}^{(3s-1)N\ell_t}}\cdot\frac{4^{\ell_t s}}{q^{2s}_{\ell_t N}(w_1^{(t)},\dots,w_{\ell_t}^{(t)})}\notag\\
		% (N\ell_t\geq n_t(1-\delta),\ s<1)\quad	& \leq \prod_{t=1}^{k-1}\frac{1}{\tilde{B}^{4sn_t(1-\delta)}} \cdot\frac{4^{\ell_t s}}{q^{2s}_{\ell_t N}(w_1^{(t)},\dots,w_{\ell_t}^{(t)})}\notag\\
		(N\ell_t\geq n_t(1-\delta),\ s<1)\quad
		& \leq \prod_{t=1}^{k-1}\biggl(\frac{1}{\tilde{B}^{4n_t}\cdot q^{2}_{\ell_t N}(w_1^{(t)},\dots,w_{\ell_t}^{(t)})}\biggr)^{s(1-\delta)}.
		\label{eq:mu-nk+2}
	\end{align}
	Note that for any $x \in E^{N}_{M}(\tilde{B})$, we have
	\begin{align}
		q_{n_k+2}(x)
		\leq 4\,\tilde{B}^{2 n_k} q_{n_k}(x)
		&\leq 4^{3N}\tilde{B}^{2 n_k}
		q_{\ell_k N}(w^{(k)}_1,\dots,w^{(k)}_{\ell_k})\, q_{n_{k-1}+2}(x)\notag\\
		&\leq \prod_{t=1}^{k}4^{3N}\tilde{B}^{2 n_t}
		q_{\ell_t N}(w^{(t)}_1,\dots,w^{(t)}_{\ell_t}).\label{eq:qnk}
	\end{align}
	It follows from Eqs.~\eqref{eq:mu-nk+2} and \eqref{eq:qnk} that 
	\begin{align}
		\mu(J_{n_{k-1}+2}(a_1,\dots,a_{n_{k-1}+2})) 
		% & \leq \prod_{t=1}^{k}\biggl(\frac{1}{\tilde{B}^{4n_t}\cdot q^{2}_{\ell_t N}(w_1^{(t)},\dots,w_{\ell_t}^{(t)})}\biggr)^{s(1-\delta)} 
		& \leq \biggl(\frac{4^{6N(k-1)}}{q^{2}_{n_{k-1} +2}}\biggr)^{s(1-\delta)}\notag \\
		& \leq \biggl(\frac{1}{q^{2}_{n_{k-1} +2}}\biggr)^{s(1-\delta)-\delta}\leq 8 |J_{n_{k-1} + 2}|^{s(1-\delta)-\delta} \label{eq:mu-J_nk-1}
	\end{align}
	where the second inequality holds since $n_k\geq \sum_{t=1}^{k-1}\ell_t N\geq 12(k-1)N^2$ by Eq.~\eqref{eq4.3} and $q_{n_{k-1}+2}^2\geq 4^{(n_{k-1}+1)/2}\geq 4^{6(k-1)N/\delta}$, and the last estimate follows from Eq.~\eqref{eqjn1}.
	\medskip
	
	(2) When $n_{k-1}+2 < n \leq m_{k}=n_{k-1} + 2 + \ell_{k} N $, 
	we first consider $n':= n_{k-1}+2 + \ell N$ for some $1 \leq \ell \leq \ell_{k}$. From the construction of $\mu$, we have
	\begin{align}
		\mu(J_{n'})
		&= \mu(J_{n_{k-1}+2})
		\prod_{i=1}^{\ell}\mu(I_N(w_i^{(k)}))\notag\\
		(\text{by Eq.~\eqref{eq:mu-J_nk-1} and Lem \ref{lem:prod-mu-i}})\quad & \leq \biggl(\frac{1}{q^2_{n_{k-1}+2}}\biggr)^{s(1-\delta)-\delta}
		\frac{1}{\tilde{B}^{(3s-1)N\ell}}\biggl(\frac{1}{q_{\ell N}^2(w^{(k)}_1,\dots,w^{(k)}_{\ell})}\biggr)^{s(1-\delta)} \notag\\
		% \frac{\tilde{B}^N}{u}\biggl(\frac{1}{\tilde{B}^{3N} q_{N}^2(w^{(k)}_i)}\biggr)^s \notag\\
		% (u\geq 1 \text{ and }3s-1>0)\quad&\leq 
		% \mu(J_{n_{k-1}+2})
		% \prod_{i=1}^{\ell}
		% \biggl(\frac{1}{q_{N}^2(w^{(k)}_i)}\biggr)^s \notag\\
		% (\text{ by Lemma }\ref{eq:2})\quad&\lesssim
		% \biggl(\frac{1}{q^2_{n_{k-1}+2}}\biggr)^{s(1-\delta)-\delta}
		% \biggl(\frac{4^\ell}{q_{\ell N}^2(w^{(k)}_1,\dots,w^{(k)}_{\ell})}\biggr)^{s} \notag\\
		% (q_{\ell N}^2\geq 4^{(\ell N-1)/2}\geq 4^{\ell/\delta})\quad&\leq
		% \biggl(\frac{1}{q^2_{n_{k-1}+2}}\biggr)^{s(1-\delta)-\delta}
		% \biggl(\frac{1}{q_{\ell N}^2(w^{(k)}_1,\dots,w^{(k)}_{\ell})}\biggr)^{s-\delta} \notag\\
		% (\tilde{B}>1 \text{ and }3s-1>0)\quad & \leq \Bigl(\frac{1}{q^2_{n_{k-1}+2}}\Bigr)^{s(1-\delta)-\delta}
		% \Bigl(\frac{1}{q_{\ell N}^2(w^{(k)}_1,\dots,w^{(k)}_{\ell})}\Bigr)^{s(1-\delta)} \notag\\
		(\text{by Lemma }\ref{eq:2}~(2))\quad&\leq \frac{1}{\tilde{B}^{(3s-1)N\ell}}\Bigl(\frac{1}{q^2_{n'}}\Bigr)^{s(1-\delta)-\delta}\label{eq:mu-lk}\\
		& \leq 8 |J_{n'}|^{s(1-\delta)-\delta}\label{eqjn'}. 
	\end{align}
	For every $n\in [n', n'+N]$ with $\ell<\ell_k$, since $a_n\leq M$, we have
	\[
	\frac{|J_n(a_1,\dots,a_n)|}{|J_{n'}(a_1,\dots,a_{n'})|}
	\geq\frac{q_{n'}^2(a_1,\dots,a_{n'})}{q_n^2(a_1,\dots,a_{n'},\dots,a_n)}
	\geq\frac{1}{16}\prod_{i=1}^{N}\frac{1}{(a_{n'+i}+1)^2}\geq \frac{1}{16}\Bigl(\frac{1}{M+1}\Bigr)^{2N}.
	\]
	Hence, by Eq.~\eqref{eqjn'},
	\[
	\mu(J_n)
	\leq \mu(J_{n'})
	\leq 8 |J_{n'}|^{s(1-\delta)-\delta}\leq 128(M+1)^{2N}|J_n|^{s(1-\delta)-\delta}\leq|J_n|^{s(1-\delta)-2\delta}.
	\]
	
	(3) When $n_{k-1} + 2 + \ell_k N=m_{k} < n \leq n_{k}=m_{k}+N+i_{k},\ i_{k}\leq N$, we have
	
	\begin{align}
		\mu(J_{n_{k}})
		&= \mu(J_{n_{k-1}+2}) 
		\prod_{i=1}^{\ell_k}\mu(I_N(w_i^{(k)})) \notag\\
		& \leq \frac{1}{\tilde{B}^{(3s-1)N\ell_k}}\Bigl(\frac{1}{q^2_{m_{k}}}\Bigr)^{s(1-\delta)-\delta} &(\text{by Eq.~\eqref{eq:mu-lk} for } \ell=\ell_k)\notag\\
		%& \leq \frac{1}{\tilde{B}^{(3s-1)N\ell_k}}\Bigl(\frac{1}{q^2_{n_{k}}}\Bigr)^{s(1-\delta)-2\delta}\notag\\
		& \leq \frac{1}{\tilde{B}^{n_k(3s-1)(1-\delta)}}\Bigl(\frac{1}{q^2_{n_{k}}}\Bigr)^{s(1-\delta)-2\delta} & (N\ell_k\geq n_k(1-\delta))\notag\\
		& \leq 3 |J_{n_{k}}|^{s(1-\delta)-2\delta}. & (\text{by Eq.}~\eqref{eqjn3})\label{J_{n_k}}
	\end{align}
	
	For $m_k<n<n_k$, we have 
	\[\mu(J_n)=\mu(J_{n_k})\leq 3|J_{n_k}|^{s(1-\delta)-2\delta}\leq 3|J_{n}|^{s(1-\delta)-2\delta}.\]
	
	(4) When $n = n_k + 1$, noticing that $q_{n_{k}+1}=a_{n_k+1}q_{n_k}+q_{n_k-1}\leq (a_{n_k+1}+1)q_{n_k}$, we have $\frac{1}{q_{n_k}}\leq \frac{a_{n_k+1}+1}{q_{n_{k}+1}}\leq \frac{2\tilde{B}^{n_k}+1}{q_{n_{k}+1}}$. Thus,
	\begin{align*}
		\mu(J_{n_{k}+1})=\frac{1}{(\# P_{n_{k}})}\cdot\mu(J_{n_{k}}) & \leq \Bigl(\frac{1}{\tilde{B}^{n_k}}\Bigr)^{1-\delta}\cdot\frac{1}{\tilde{B}^{n_k(3s-1)(1-\delta)}}\Bigl(\frac{1}{q^2_{n_{k}}}\Bigr)^{s(1-\delta)-2\delta}\\
		& \leq \Bigl(\frac{1}{\tilde{B}^{3n_k}}\Bigr)^{s(1-\delta)}\Bigl(\frac{(2\tilde{B}^{n_k}+1)^2}{q^2_{n_k+1}}\Bigr)^{s(1-\delta)-2\delta}\\
		& \leq 3 |J_{n_k+1}|^{s(1-\delta)-2\delta}
	\end{align*}
	where the first inequality follows from Proposition \ref{p2.5} and Eq.~\eqref{J_{n_k}}, and the last step follows from Eq.~\eqref{eqjn3}.
	
	As a consequence, we have verified that for all $n\geq n_1$, 
	\begin{equation}
		\mu(J_n) \leq 8 |J_n|^{s(1-\delta)-2\delta}. \label{eq:sec:4.2.1}
	\end{equation}
	
	\subsubsection{H\"older exponent for balls}
	We now estimate the H\"{o}lder exponents for small balls. 
	For any ball \( B(x,r) \) centered at $x \in E^{N}_{M}(\tilde{B})$, there exists a unique integer $n$ such that
	$$|G_{n+1}(a_1,\dots,a_{n+1})|\leq r<|G_{n}(a_1,\dots,a_{n})|$$
	where $x$ has the continued fraction expansion $(a_1,a_2,\dots)$, $\mathbf{a}:=(a_1,\dots,a_n)$ and $G_n(\mathbf{a})$ is the minimal gap between $J_n(\mathbf{a})$ and its adjacent fundamental intervals $J_{n}(\mathbf{a}')$ with $\mathbf{a}'\in D_n$ and $\mathbf{a}'\neq \mathbf{a}$. A lower bound for the length of $G_n$ is given below.
	\begin{pro}[\cite{TTW}]\label{prop:gap}
		With the notation above, we have
		\[
		|G_n(a_1,\dots,a_{n})|
		\ge \frac{1}{8M} \, |J_n(a_1,\dots,a_{n})|.
		\]
	\end{pro}

	We now evaluate $\mu(B(x,r))$. Since $r<|G_{n}(a_1,\dots,a_{n})|$, \( B(x,r) \) intersects only one $J_n$, that is $J_n(\mathbf{a})$. We always have $\mu(B(x,r))\leq \mu(J_n(\mathbf{a}))$.
	
	\textbf{Case 1:} When $n_{k-1}+2\leq n<n_k-1$, we see that $a_{n+1}\leq M$.  By Eqs.~\eqref{eqjn1} and \eqref{eqjn2}, we have 
	\begin{equation}\label{eq:jn-jn+1}
		|J_n(\mathbf{a})|\leq \frac{2}{q_n^{2}(\mathbf{a})}\leq \frac{2(a_{n+1}+1)^2}{q_{n+1}^2(\mathbf{a},a_{n+1})}\leq 36(a_{n+1}+1)^2|J_{n+1}(\mathbf{a},a_{n+1})|.
	\end{equation}
	It follows from Eqs.~\eqref{eq:sec:4.2.1}, \eqref{eq:jn-jn+1} and Proposition~ \ref{prop:gap} that 
	\begin{align*}
		\mu(B(x,r))\leq \mu(J_n(\mathbf{a})) & \leq 8 |J_n(\mathbf{a})|^{s(1-\delta)-2\delta}\\
		& \leq 8\cdot 36\cdot (M+1)^2|J_{n+1}(\mathbf{a},a_{n+1})|^{s(1-\delta)-2\delta}\\
		&\leq 8^2\cdot 36\cdot (M+1)^3 |G_{n+1}(\mathbf{a},a_{n+1})|^{s(1-\delta)-2\delta}\leq r^{s(1-\delta)-3\delta}.
	\end{align*}

	\textbf{Case 2:} when $n=n_k-1,n_k,\,n_k+1$. 
	
	Combining Lemma \ref{L2.9} and the fact that $a_{n+1}q_{n}\leq q_{n+1}\leq (a_{n+1}+1)q_{n}$, we see that for any $x\in I_{n+1}(a_1,\dots,a_{n+1})$ with $a_n\geq 2$, $$B\bigl(x,|I_{n+1}(a_1,\dots,a_{n+1})|\bigr)\subset \bigcup_{i=-1}^{2}I_{n+1}(a_1,\dots,a_{n},a_{n+1}+i).$$ Write $I_{n+1}^{(i)}:=I_{n+1}(a_1,\dots,a_{n},a_{n+1}+i)$ and define $J_{n+1}^{(i)}$ similarly.. 
	
	If $r\leq |I_{n+1}^{(0)}|$, the ball $B(x,r)$ intersects at most 4 level-$(n+1)$ cylinders $I_{n+1}^{(-1)},\dots,I_{n+1}^{(2)}$. Note that by Eq.~\eqref{eqjn3}, we have $|J_{n+1}^{(i)}|\leq 2|J_{n+1}^{(0)}|$. Then, we have
	\begin{align*}
		\mu(B(x,r))&\leq \mu(B(x,|I_{n+1}^{(0)}|))\leq \sum_{i=-1}^{2}\mu(J_{n+1}^{(i)})\leq \sum_{i=-1}^{2} 8 |J_{n+1}^{(i)}|^{s(1-\delta)-2\delta}\\
		& \leq 64|J_{n+1}^{(0)}|^{s(1-\delta)-2\delta} \leq 512 M|G_{n+1}|^{s(1-\delta)-2\delta}\leq r^{s(1-\delta)-3\delta}.
	\end{align*}
	
	% {\color{blue}
		% We next estimate the number of level-$(n+1)$ basic cylinders intersecting $B(x,r)$. Let
		% \[
		% \Lambda_{+i}:=|\cup_{k=0}^{i}I_{n+1}(a_1,\dots,a_{n},a_{n+1}+k)|.
		% \]
		% By Lemma \ref{L2.9}, we have $\Lambda_{+i}=\frac{i+1}{q_{n+1}(q_{n+1}+(i+1)q_{n})}$.
		% Setting $\Lambda_{+i}=r$, we obtain $i=\frac{rq_n^2}{1-rq_nq_{n-1}}-1$.
		% If	$r\geq |I_n(a_1,\dots,a_n)|=\frac{1}{q_n(q_n+q_{n-1})}$, it follows that $rq_nq_{n-1}\geq \frac{q_{n-1}}{q_n+q_{n-1}}
		% \geq \frac{1}{a_n+2}.$
		% The same estimate clearly holds for $\Lambda_{-i}$. Hence, if $a_n$ is bounded, then there exists a constant $C>0$, the ball $B(x,r)$ intersects at most $C\,\frac{r}{|I_n(a_1,\dots,a_n)|}$ 
		% basic cylinders of level-$n$.	}
	
	Now suppose $r> |I_{n+1}^{(0)}|$. For $n=n_k$, by the construction of $E_M^{N}(\tilde{B})$, the length of cylinders $I_{n+1}^{(i)}$ in $E_M^{N}(\tilde{B})$ are comparable, since $a_{n+1}+i\in [\tilde{B}^{n_k},2\tilde{B}^{n_k}]$. Thus, the ball $B(x,r)$ intersects at most $C\frac{r}{|I_{n+1}^{(0)}|}$ cylinders $I_{n+1}^{(i)}$. Since $|I_{n+1}^{(0)}|\geq \frac{1}{8\tilde{B}^{2n_k}q^2_{n_k}}$, we have  
	\begin{align*}
		\mu(B(x,r))
		&\leq \min\bigl\{\mu(J_{n}),(8Cr\tilde{B}^{2n_k}q^2_{n_k}) \mu(J_{n+1})\bigr\}\\
		&\leq \mu(J_{n})\min\bigl\{1,(8Cr\tilde{B}^{2n_k}q^2_{n})\tfrac{1}{\# P_{n_k}}\bigr\} & (\text{by the definition of }\mu)\\
		&\leq \mu(J_{n})\min\{1,\ 8Cr\tilde{B}^{n_k(1+\delta)}q^2_{n_k}\}& \text{(by Proposition \ref{p2.5})}\\
		&\leq 8\Bigl(\frac{1}{\tilde{B}^{n_k}q^2_{n_k}}\Bigr)^{s(1-\delta)-2\delta}(8Cr\tilde{B}^{n_k(1+\delta)}q^2_{n_k})^{s(1-\delta)-2\delta} & \text{(by Eqs.~\eqref{eq:sec:4.2.1} \& \eqref{eqjn3})}\\
		&\leq r^{s(1-\delta)-3\delta}.
	\end{align*}
	The case $n=n_k+1$ is similar since $\mu(J_{n+1})=\frac{1}{\#P_{n_k}}\mu(J_n)$. For $n=n_k-1$, we have $\mu(J_{n+1})=\mu(J_n)$ and the same argument applies.
	
	% Letting $b$ be the prime next to $a_{n_k+1}$, by Eq. \eqref{eqjn3}, we have 
	% \[\frac{1}{3}\frac{q^2_{n_k+1}(\mathbf{a},b)}{q^2_{n_k+1}(\mathbf{a},a_{n_k+1})}\leq \frac{|J_{n_k+1}(\mathbf{a},a_{n_k+1})|}{|J_{n_k+1}(\mathbf{a},b)|} \leq  3\frac{q^2_{n_k+1}(\mathbf{a},b)}{q^2_{n_k+1}(\mathbf{a},a_{n_k+1})}\]
	% where $\mathbf{a}=(a_1,\dots,a_{n_k})$. Since, by Proposition \ref{p2.5}~(2), $a_{n_k+1} \geq 0.999b$ or $b\geq 0.999 a_{n_k+1}$,
	% \[\frac{q^2_{n_k+1}(\mathbf{a},b)}{q^2_{n_k+1}(\mathbf{a},a_{n_k+1})}= \Bigl(\frac{bq_{n_k}+q_{n_k-1}}{a_{n_k+1}q_{n_k}+q_{n_k-1}}\Bigr)^2\in \Bigl[\frac{b}{a_{n_k+1}+1},\ \frac{b+1}{a_{n_k+1}}\Bigr]\subseteq [\tfrac{1}{2},\ 2].\]
	% Consequently, \[\tfrac{1}{6}|J_{n_k+1}(\mathbf{a},b)|\leq |J_{n_k+1}(\mathbf{a},a_{n_k+1})|\leq 6 |J_{n_k+1}(\mathbf{a},b)|.\]
	
	Therefore, for all $x \in E^{N}_{M}(\tilde{B})$, for all small $r > 0$, $\mu(B(x,r)) \leq  r^{s(1-\delta)-3\delta}$.
	Applying the Mass Distribution Principle \cite[Page 67]{Falconer}, we conclude that for all $1/2 < s < s(B)$,
	\[
	\dim_H E'(\phi) \geq \dim_H E^N_M(\tilde{B}) \geq s(1-\delta) - 3\delta.
	\]
	Letting $\delta \to 0$, and using the arbitrariness of $s$, we obtain
	$\dim_H E'(\phi) \geq s(B)$.

	\section{Proof of Theorem \ref{t2}: the other case}
	In this section, we complete the proof of Theorem~\ref{t2}. Recall that $\phi_{b,c}(n)=c^{b^n}$ and define
	\begin{align*}
		&F'(\phi_{b,c}) := \{x \in [0,1) : a'_n(x) \geq c^{b^n} \ \text{for i.m. } n \in \mathbb{N}\},\\
		&F''(\phi_{b,c}) := \{x \in [0,1) : a'_n(x) \geq c^{b^n} \ \text{for all } n \in \mathbb{N}\}.
	\end{align*}
	It can be readily verified that these sets satisfy the following inclusions:
	\[
	F''(\phi_{b,c})
	\subset
	\{x \in [0,1) : a'_{n-1}(x) \geq c^{b^n}, \, a'_n(x) \geq c^{b^n} \ \text{for i.m. } n \in \mathbb{N}\}
	\subset E'(\phi_{b,c})\subset
	F'(\phi_{b,c}).
	\]
	It is known from \cite[Theorem~1.11]{Robert} that $\dim_H F'(\phi_{b,c})=\dim_H F''(\phi_{b,c})=\frac{1}{1+b}$. Thus, we have
	\begin{equation}\label{e5.1}
		\dim_H E'(\phi_{b,c}) = \frac{1}{1+b}.
	\end{equation}
	
	\paragraph{\textbf{Case 1: \( B = 1 \).}}
	For any \( \sigma\in (1,+\infty) \), for $n\in \mathbb{N}$, define $\phi_\sigma(n) = \sigma^n \phi(n)$.
	Then
	\[
	\log B = \liminf_{n \to \infty} \frac{\log \phi_\sigma(n)}{n} = \log \sigma > 0,
	\qquad
	E'(\phi_\sigma) \subset E'(\phi).
	\]
	Now $B=\sigma$. By Theorem \ref{t2} (the subcase $1<B<+\infty$) proved in Section 4, we have $\dim_H E'(\phi_\sigma) \geq s(B) = s(\sigma)$. Since $E'(\phi)\supset E'(\phi_{\sigma})$, we see that
	\[
	\dim_H E'(\phi)\geq\dim_H E'(\phi_\sigma) \geq s(\sigma).
	\]
	By Lemma \ref{pro4.1}(4), letting \( \sigma \to 1 \) yields \( \dim_H E'(\phi) = 1 \).
	\medskip
	
	\paragraph{\textbf{Case 2: \( B = \infty \).}}
	In this case, we distinguish the following two subcases.
	
	\subparagraph{Case 2.1: \( b = 1 \).}
	For any \( B > 1 \), set \( \phi_B(n) = B^n \) for all \( n \geq 1 \).
	On the one hand, it is evident that
	\[
	E'(\phi) \subset E'(\phi_B), \qquad \text{for all } 1 < B < \infty.
	\]
	Hence, we have $\dim_H E'(\phi) \leq \tfrac{1}{2}$ by Lemma~\ref{pro4.1}.
	%	\[
	%	\dim_H E'(\phi) \leq \tfrac{1}{2}.
	%	\]
	
	On the other hand, since $b=1$, for any \( \lambda > 1 \), we have
	\[
	\phi(n) \leq e^{\lambda^n} \quad \text{for i.m. } n \in \mathbb{N}.
	\]
	Consequently,
	\[
	\{x \in [0,1) : a'_n(x) \geq e^{\lambda^{n}} \ \text{for all } n \in \mathbb{N}\}
	\subset E'(\phi).
	\]
	Therefore, by~\eqref{e5.1},
	\[
	\dim_H E'(\phi) \geq \frac{1}{\lambda + 1}.
	\]
	Letting \( \lambda \to 1 \) gives the desired result.
	
	% \subparagraph{Case 2.2: \( b = \infty \).}
	% For any finite \( \lambda \), we have
	% \[
	% E'(\phi) \subset
	% \Bigl\{ x \in [0,1) : \exists\, 1 \leq k < n, \, a'_k(x) \geq e^{\lambda^n}, \, a'_n(x) \geq e^{\lambda^n} \ \text{for i.m. } n \in \mathbb{N} \Bigr\}.
	% \]
	% Therefore, by~\eqref{e5.1}, we obtain
	% \[
	% \dim_H E'(\phi) = 0.
	% \]
	\medskip
	\subparagraph{Case 2.2: \( 1 < b \leq \infty \).}
	For any \( \lambda < b \), we have $\phi(n) \geq e^{\lambda^n}$ for all sufficiently large $n$. Thus,
	\[
	E'(\phi) \subset
	\Bigl\{ x \in [0,1) : \exists\, 1 \leq k < n, \, a'_k(x) \geq e^{\lambda^n}, \, a'_n(x) \geq e^{\lambda^n} \ \text{for i.m. } n \in \mathbb{N} \Bigr\}=E'(\phi_{\lambda,e}),
	\]
	and by~\eqref{e5.1}, we have  \( \dim_H E'(\phi_{\lambda,e})=\frac{1}{\lambda + 1} \). Then $\dim_H E'(\phi)\leq \dim_H E'(\phi_{\lambda,e})=\frac{1}{\lambda + 1}$. Letting $\lambda \to b$, we have $\dim_H E'(\phi)\leq \frac{1}{b+1}$.
	
	When $b=\infty$, letting $\lambda \to \infty$, we have $\dim_H E'(\phi)=0$.  When $b<\infty$, for any \( \theta > b \), one has
	\[
	\phi(n) \leq e^{\theta^n} \quad \text{for i.m. } n \in \mathbb{N}.
	\]
	It follows that
	\[
	F''(\phi_{\theta, e}) = \{ x \in [0,1) : a'_n(x) \geq e^{\theta^{n+1}} \ \text{for all } n \in \mathbb{N} \}
	\subset E'(\phi),
	\]
	and $\dim_H E'(\phi)\geq \dim_H F''(\phi_{\theta, e}) = \frac{1}{\theta + 1}$. Letting $\theta \to b$, we have $\dim_H E'(\phi)\geq \frac{1}{b+1}$. This ends the proof.
	
	\medskip
	\noindent\textbf{Acknowledgment.} This work is supported by the National Natural Science Foundation of China (No.12371086, 12271175) and the National Key R\&D Program of China (No. 2024YFA1013701).
		
\end{document}